\documentclass[11pt]{article}
\usepackage[margin=1.03in]{geometry}
\usepackage{amsmath,amssymb,amsthm,mathtools}
\usepackage[bookmarksnumbered,hidelinks]{hyperref}
\hypersetup{pdftitle={Embeddable matrix-algebra bundles, tensor units and UHF stabilization},pdfauthor={A. V. Ershov}}
\newtheorem{theorem}{Theorem}[section]
\newtheorem{proposition}[theorem]{Proposition}
\newtheorem{lemma}[theorem]{Lemma}

\theoremstyle{definition}
\newtheorem{definition}[theorem]{Definition}
\newtheorem{example}[theorem]{Example}
\theoremstyle{remark}
\newtheorem{remark}[theorem]{Remark}
\newcommand{\Z}{\mathbb Z}
\newcommand{\C}{\mathbb C}
\newcommand{\Q}{\mathbb Q}
\newcommand{\RP}{\mathbb {RP}}
\newcommand{\CP}{\mathbb {CP}}
\newcommand{\Gr}{\mathrm{Gr}}
\newcommand{\Fr}{\mathrm{Fr}}
\newcommand{\PU}{\mathrm{PU}}
\newcommand{\SU}{\mathrm{SU}}
\newcommand{\BSU}{B\mathrm{SU}}
\newcommand{\BU}{B\mathrm U}
\newcommand{\UU}{\mathrm U}
\newcommand{\Aut}{\mathrm{Aut}}
\newcommand{\End}{\mathrm{End}}
\newcommand{\Emb}{\mathrm{Emb}_1}
\newcommand{\Tel}{\operatorname{Tel}}
\newcommand{\hofib}{\operatorname{hofib}}
\newcommand{\rk}{\operatorname{rk}}
\newcommand{\id}{\mathrm{id}}
\newcommand{\one}{\mathbf{1}}
\newcommand{\Bun}{\mathcal B}
\newcommand{\FAB}{\mathrm{FAB}_{\mathrm{st}}}
\newcommand{\units}{\mathcal U_{\mathrm{SU}}}
\newcommand{\Kc}{\mathcal K}
\newcommand{\arx}[1]{\href{https://arxiv.org/abs/#1}{\texttt{arXiv:#1}}}
\allowdisplaybreaks[2]
\title{Embeddable matrix-algebra bundles,\\ tensor units and UHF stabilization}
\author{A.\,V. Ershov}
\date{}
\begin{document}
\maketitle
\begin{abstract}
A unital embedding of a bundle of complex matrix algebras in a trivial matrix algebra bundle of coprime multiplicity determines two complementary factors. We recover from them an integral virtual special-unitary bundle of rank one. This gives a proof of the stable classification by $\BSU_{\otimes}$ that allows torsion in the cohomology of the base. At fixed core rank $k$, the telescope of matrix Grassmannians with complementary ranks $\ell^N$ is the homotopy fibre of $B\PU(k)\to B\Aut(M_k\otimes M_{\ell^\infty})$. Thus unital UHF triviality detects exactly which cores admit coprime finite-dimensional embeddings. We determine the least stage on even spheres and on four-dimensional Moore spaces. A rank-three example over $S^6$ first embeds in $M_{12}$ and represents a unit whose clutching gives a higher twist over $S^7$.
\end{abstract}

\section{Introduction}

A family of subalgebras $M_k\subset M_{km}$ determines an $M_k$-bundle and a complementary $M_m$-bundle. Multiplication identifies their tensor product with the trivial $M_{km}$-bundle. When $k$ and $m$ are coprime, these two pieces of data recover an integral class in complex $K$-theory. The class records how the matrix factors vary even though the ambient algebra is trivial. This construction was introduced under the name \emph{floating algebra bundle} (FAB) in \cite{short,old}. We use \emph{embeddable} for the property of a core and \emph{embedded pair} for a core together with its embedding.

Two stabilization procedures answer different questions. When both ranks increase, stable equivalence classes form the group of virtual rank-one special-unitary bundles under tensor product. We give a proof using complementary factors in Sections~\ref{sec:geometry}--\ref{sec:classification}. At fixed $k$ one can instead increase only the complementary rank. For any $\ell\ge2$ prime to $k$, this gives a criterion for a single $M_k$-bundle $\Bun$ over a finite CW complex $X$:
\[
 \Bun\text{ admits a coprime embedding}
 \quad\Longleftrightarrow\quad
 \Bun\otimes M_{\ell^\infty}\cong X\times(M_k\otimes M_{\ell^\infty}).
\]
Section~\ref{sec:uhf} identifies the homotopy fibre that gives this criterion. The argument combines finite matrix frames, already studied in \cite{brauer,obstruct}, with the contractibility theorem of Dadarlat and Pennig \cite{DP}. It uses ordinary bundles and mapping telescopes.

The distinction between an embedding and its abstract core is essential on bases with torsion. Corollary~22 of \cite{old} asserted uniqueness of a special-unitary lift of the core. Example~\ref{ex:RP} gives two lifts of the trivial algebra bundle with different second Chern classes. The proof here constructs a class from both factors, checks its invariance under automorphisms of the ambient trivial algebra, and defines characteristic classes before pulling back to the base.

Sections~\ref{sec:stages} and~\ref{sec:twists} consider the size of a finite realization. For a rank-$k$ vector bundle of nonzero Bott integer $t$ on $S^{2j}$, $2\le j\le k$, the multiplicities that work are exactly the $m\ge j$ with $k\mid mt$. Over torsion bases, low-rank representatives also matter. We describe all stable classes over a fixed embeddable core as a coset of a subgroup determined by $k$-torsion in $K^0(X)$. On a Moore space this gives both the number of such classes and, together with the multiplicity equation, the first embedding stage. The final example uses the suspension construction of Pennig \cite{Pennig} to relate a least matrix size to a twisted $K$-group.

\paragraph{Conventions.}
Bases are finite CW complexes unless specified. For topological groups we use numerable classifying bundles. Arguments apply separately on each component. All matrix algebras are over $\C$, all algebra embeddings are unital $*$-homomorphisms, and automorphism groups have the point-norm topology. The complex algebra convention in \cite{old} gives the same classifying homotopy types: the inclusions of the projective unitary groups in the complex projective linear groups are homotopy equivalences and induce equivalences of the corresponding homogeneous spaces. Write $M_s$ for the $s$ by $s$ matrices, $\one^s$ for a trivial rank-$s$ vector bundle and $\underline M_s$ for a trivial matrix algebra bundle. Integers inside $K^0(X)$ denote trivial bundle classes. We use the standard complex stable range: in dimension at most $2s$, every stable class of rank $s$ has a unique rank-$s$ representative up to isomorphism. This follows from the $(2s+1)$-connectivity of $\BU(s)\to\BU$, using $\UU(s+1)/\UU(s)=S^{2s+1}$. The corresponding special-unitary assertion follows by taking homotopy fibres of the determinant maps. We use Bott periodicity and the usual finite CW filtration of $K^0$ as in \cite{Atiyah,Bott}.

\section{Embedded pairs and compatible lifts}\label{sec:geometry}

An embedded pair of type $(k,m)$ over $X$ consists of an $M_k$-bundle $\Bun$ and an embedding $\Bun\to\underline M_{km}$. Its commutant $\mathcal C$ is a locally trivial $M_m$-bundle. Two pairs of the same type are isomorphic if an isomorphism of their cores and a bundle automorphism of the ambient trivial algebra intertwine the embeddings. We work with $k,m>1$ in the stable theory and also allow $m=1$ when considering least multiplicities.

Put
\[
 \Fr_{k,m}=\Emb(M_k,M_{km}),\qquad
 \Gr_{k,m}=\Fr_{k,m}/\PU(k).
\]
Unitary conjugacy of matrix units identifies $\Fr_{k,m}$ with
$\UU(km)/(1_k\otimes\UU(m))$ and $\Gr_{k,m}$ with
$\PU(km)/(\PU(k)\otimes\PU(m))$. The source action on frames is free. Its quotient is the space of the subalgebras themselves. A family in the fixed trivial algebra is a continuous map to $\Gr_{k,m}$. The bundle $\Fr_{k,m}\to\Gr_{k,m}$ is the principal frame bundle of the universal core. In particular its classifying map gives the homotopy fibre sequence
\begin{equation}\label{eq:finite}
 \Fr_{k,m}\longrightarrow \Gr_{k,m}\longrightarrow B\PU(k).
\end{equation}
Multiplication identifies the universal core tensored with its commutant with $\underline M_{km}$, by the matrix double centralizer theorem in local coordinates.

\begin{lemma}\label{lem:pair}
Suppose $(k,m)=1$. The Grassmannian has the model
\[
 \Gr_{k,m}=\SU(km)/(\SU(k)\otimes\SU(m)).
\]
There are special-unitary bundles $V,W$ on it of ranks $k,m$, equipped with a trivialization $V\otimes W\cong\one^{km}$. Their endomorphism algebras are the universal core and its commutant.
\end{lemma}
\begin{proof}
The kernel of tensor multiplication on $\SU(k)\times\SU(m)$ consists of $(z1_k,z^{-1}1_m)$ with $z^k=z^m=1$, so it is trivial. Multiplication of the two centres maps onto the centre of $\SU(km)$ by coprimality. Passing to projective groups therefore induces the stated identification of quotients. Associate to the homogeneous principal $\SU(k)\times\SU(m)$-bundle its two defining representations. Their tensor product extends to the defining representation of $\SU(km)$, giving the claimed trivialization. The associated endomorphism bundles act on the first and second factors.
\end{proof}

\begin{example}[The special-unitary lift of a core]\label{ex:RP}
On $\RP^4$ let $L$ be the complexified tautological real line and $a=c_1(L)$. We have $L^2\cong\one$ (powers of a line denote tensor powers) and $a^2\ne0$ in $H^4(\RP^4;\Z)=\Z/2$. Indeed the mod-two reduction of $a$ is the square of the degree-one generator, by the Stiefel--Whitney formula for the underlying real plane of $L$. Its square reduces to the nonzero fourth power. Both $V_0=\one^2$ and $V_1=L\oplus L$ have trivial determinant and trivial endomorphism algebra, whereas $c_2(V_0)=0$ and $c_2(V_1)=a^2$. Their common core has the standard embedding in $\underline M_6$. There is no $\SU(3)$-bundle $W'$ with $V_1^{\oplus3}\cong(W')^{\oplus2}$, since the second Chern classes would be $a^2$ and $0$. Thus the existence assertion in Lemma~21 of \cite{old} holds, while the uniqueness assertion of its Corollary~22 and the arbitrary-lift formulation of its Remark~32 fail.
\end{example}

Let $R_X=K^0(X)$ and set
\[
 \units(X)=\{u\in R_X:\rk u=1,\ c_1(u)=0\}.
\]
The rank-zero ideal is nilpotent by the skeletal filtration. Hence $\units(X)$ is a group under multiplication and is represented by $\BSU_\otimes$, the usual $\BSU$ with the tensor product of virtual rank-one classes. The determinant splitting for stable complex bundles identifies $[X,\BSU]$ with the rank-zero classes of first Chern class zero, translated here by $1$.

\begin{lemma}[Integral normalization]\label{lem:normalize}
Let $V,W$ be as in Lemma~\ref{lem:pair}, or their pullbacks. There is exactly one $u\in\units(X)$ such that
\begin{equation}\label{eq:two}
 [V]=ku,\qquad [W]=mu^{-1}.
\end{equation}
The class is natural in $X$ and invariant under isomorphisms of embedded pairs.
\end{lemma}
\begin{proof}
For any abelian group $P$ and coprime $k,m$, the diagonal identifies $P$ with the pullback of $P[1/k]$ and $P[1/m]$ over $P[1/(km)]$. For injectivity use B\'ezout on the powers that kill a given element. For surjectivity represent a matching pair by $a/k^i,b/m^j$. After increasing the powers, $m^j a=k^i b$ in $P$. If $h k^i+t m^j=1$, then $h a+t b$ has the prescribed images.

Apply this to $R_X$. The classes $[V]/k$ in $R_X[1/k]$ and $([W]/m)^{-1}$ in $R_X[1/m]$ agree after inverting $km$ because $[V][W]=km$. Each has rank one and hence is a unit. They give a unique rank-one unit $u$. The two identities \eqref{eq:two} hold in both localizations by the same product equation, and hence hold integrally. They imply $k c_1(u)=m c_1(u)=0$, which gives $c_1(u)=0$. They also force the two localized values and prove uniqueness.

Projectively isomorphic vector bundles of equal rank differ by a tensor factor which is a line bundle: local unitary intertwiners differ by scalars and their spans glue into the line of intertwiners. Thus, for two choices of first factor with trivial determinants and isomorphic matrix algebras, the line $L$ satisfies $L^k=\one$. In $R_X[1/k]$ this implies $[L]=1$: the factor $1+[L]+\cdots+[L]^{k-1}$ has rank $k$ and is a unit there. The first localized value is consequently independent of this ambiguity. The complementary projective factor similarly determines the second localized value. An isomorphism of embedded pairs preserves both projective factors, and pullback respects the construction.
\end{proof}

Swapping the factors takes $u$ to $u^{-1}$. Tensoring embedded pairs takes $(u,u')$ to $uu'$ whenever the product ranks are admissible: the two tensor products of vector factors satisfy \eqref{eq:two}, and uniqueness applies. In particular, amplification by a trivial embedded pair preserves $u$.

\section{Stable classification and characteristic classes}\label{sec:classification}

For $k\mid k'$ and $m\mid m'$ the amplification $\Gr_{k,m}\to\Gr_{k',m'}$ enlarges the core by a trivial $M_{k'/k}$ and the complement by a trivial $M_{m'/m}$. Tensor coordinates are permuted to keep the two factors in this order. For a coprime pair let $q_{k,m}:\Gr_{k,m}\to\BSU$ represent its universal class $u$, in rank-one coordinates. The matrix Grassmannian is a compact smooth manifold and hence has finite CW type.

\begin{proposition}\label{prop:stable-range}
If $(k,m)=1$, the map $q_{k,m}$ induces isomorphisms on $\pi_i$ for $i\le 2\min(k,m)$. For $\dim X\le 2\min(k,m)$ the invariant $u$ classifies embedded pairs of type $(k,m)$ up to isomorphism, and every class in $\units(X)$ occurs.
\end{proposition}
\begin{proof}
Put $b=\min(k,m)$. In the homotopy sequence of $\SU(k)\times\SU(m)\to\SU(km)\to\Gr_{k,m}$, the map in degree $2j-1$, $2\le j\le b$, is
\[
 \Z\oplus\Z\longrightarrow\Z,\qquad (x,y)\longmapsto mx+ky .
\]
On each factor the representation is the direct sum of $m$ or $k$ copies, which gives the asserted multiplication under stable clutching. It is surjective. Bott periodicity gives $\pi_{2j}\SU(km)=0$, and the even groups of both factors in degree $2j-2$ vanish in the indicated range. Hence $\pi_{2j-1}\Gr_{k,m}=0$ and the boundary identifies $\pi_{2j}\Gr_{k,m}$ with the kernel generated by $(k,-m)$. In degrees one and two the homotopy groups vanish by simple connectivity and $\pi_2\SU(s)=0$. Thus the space is $3$-connected. Over a generating $2j$-sphere the reduced class of the first vector factor is $k$ times a Bott generator, up to a common sign. Equation~\eqref{eq:two} now says that $q_{k,m}$ takes the generator to a generator of $\pi_{2j}\BSU$.

A map inducing isomorphisms through degree $d$ between these simply connected CW types induces an equivalence of their $d$th Postnikov sections, hence a bijection on homotopy classes from any CW complex of dimension at most $d$. Apply this with $d=2b$. Every embedded family in the fixed trivial algebra determines a map to the Grassmannian. Homotopic maps give isomorphic embedded pairs: the principal bundle of adapted frames over $X\times I$ has isomorphic restrictions at the ends, and such an isomorphism extends to the associated ambient matrix bundles. Conversely isomorphic pairs have the same $u$ by Lemma~\ref{lem:normalize}, so their maps are homotopic in this range.
\end{proof}

For comparison, for arbitrary $k,m>1$ put $g=(k,m)$. The projective group fibration gives $\pi_1\Gr_{k,m}=\pi_2\Gr_{k,m}=\Z/g$. Indeed these are the cokernel and kernel of $(x,y)\mapsto mx+ky$ from $\Z/k\oplus\Z/m$ to $\Z/(km)$. In the range $2\le j\le\min(k,m)$ one obtains $\pi_{2j}\Gr_{k,m}=\Z$ and $\pi_{2j-1}\Gr_{k,m}=\Z/g$ from the corresponding map on $\Z^2$. For $(k',m')=(ak,bm)$ write $g'=(k',m')$. The even-degree generator of the kernel $(k/g,-m/g)$ maps to $(k'/g,-m'/g)$, giving multiplication by $g'/g$ in all the indicated even degrees, with residue classes in degree two. The odd-degree map sends the class of $1$ to $ab\bmod g'$ by ambient amplification. This recovers the finite homotopy group computations of \cite{old}.

\begin{definition}
The product of pairs is formed by tensoring their cores, their complements and their ambient algebras. Let $I_{k,m}$ denote the constant pair with core $M_k\otimes1_m$. We impose the stable equivalence relation of \cite[Definition~15]{old}: it is generated by identifying pairs $P$ of type $(k,m)$ and $P'$ of type $(k',m')$ when $(kk',mm')=1$ and $P\otimes I_{k',m'}$ is isomorphic to $P'\otimes I_{k,m}$. Denote the resulting set of classes by $\FAB(X)$.
\end{definition}

\begin{theorem}\label{thm:stable}
On finite CW complexes the invariant in \eqref{eq:two} induces a natural group isomorphism
\[
                 \FAB(X)\ \cong\ \units(X)=[X,\BSU_\otimes].
\]
For any coordinatewise divisibility tower of coprime pairs $(k_i,m_i)$ with $k_i,m_i\to\infty$, its Grassmannian telescope is homotopy equivalent to $\BSU$. With the stable tensor product it represents $\BSU_\otimes$. Centralizers represent inversion.
\end{theorem}
\begin{proof}
The invariant respects the generating equivalences. Any given class in $\units(X)$ is represented in every coprime pair of sufficiently large ranks by Proposition~\ref{prop:stable-range}. To compare two embedded pairs having this class, choose a third pair of large distinct prime ranks $(r,s)$ avoiding every prime in the four given ranks. Realize the same class in type $(r,s)$. After cross-amplifying with this third type, the resulting two maps into $\Gr_{kr,ms}$ have equal invariant and are in the dimension range of Proposition~\ref{prop:stable-range}. Their embedded pairs are isomorphic. Do this separately for each original pair. The third representative then connects them by the prescribed relation. This proves injectivity even when the two original types cannot be amplified to a common coprime type directly. Taking admissible representatives, the tensor and centralizer formulas following Lemma~\ref{lem:normalize} give the group structure and inversion.

The $q_{k_i,m_i}$ are compatible up to homotopy because the normalized classes agree on each stage. Choosing homotopies on successive mapping cylinders defines a map from the telescope. It induces isomorphisms on every homotopy group by Proposition~\ref{prop:stable-range}. The spaces have CW type, so this is a homotopy equivalence. Transport of tensor multiplication along it gives an $H$-structure. On every pair of admissible finite stages its product agrees up to homotopy with tensoring the two families: on their product the two maps to $\BSU$ represent the same class by \eqref{eq:two}. For a fixed tower, a product type need not belong to the tower. Its universal class still gives the comparison to $\BSU$, and its restriction to any finite base represents the group law just proved.
\end{proof}

For example, $(2^i,3^i)$ is a single integral model. The admissible pairs under divisibility form no directed set, since $(2,3)$ and $(3,2)$ have no common admissible multiple. The common homotopy type of the towers is the meaning of the stable Grassmannian in this theorem. For $u=1+x,\ v=1+y$, tensor product induces $x+y+xy$ on rank-zero special-unitary classes.

\subsection*{A prescribed ambient special-unitary bundle}
Let $T$ have rank $km$, with $c_1(T)=0$, $(k,m)=1$, and $\dim X\le2\min(k,m)$. Consider isomorphism classes of ordered special-unitary vector factors of ranks $k,m$ whose tensor product is isomorphic to $T$. An isomorphism with $T$ is not part of the data.

\begin{proposition}\label{prop:torsor}
Such factor pairs exist and form a torsor for $\units(X)$. On stable classes the action is $u\cdot(v,w)=(uv,u^{-1}w)$.
\end{proposition}
\begin{proof}
Let $J$ be the ideal of rank-zero classes with first Chern class zero. It is nilpotent. Start with $v=k,w=m$ and error $e=[T]-vw\in J$. Whenever $e\in J^i$, take integers $a,b$ with $ma+kb=1$ and replace $v,w$ by $v+ae,w+be$. The new error belongs to $J^{i+1}$, since $v-k,w-m\in J$. Iteration yields $vw=[T]$. The stable range realizes $v,w$ by bundles, and their tensor product is isomorphic to $T$ by cancellation.

For two solutions $(v,w),(v',w')$, the quotient $v'/v$ in $R_X[1/k]$ and the quotient $w/w'$ in $R_X[1/m]$ agree in $R_X[1/(km)]$ by the product equation. The denominators are units there by their ranks and nilpotence. They glue to a unique $u$ of rank one. The equations $v'=uv$ and $w'=u^{-1}w$ hold in both localizations and hence integrally. The two equations give $k c_1(u)=m c_1(u)=0$. This proves simple transitivity. A chosen factor pair therefore identifies the torsor with the group.
\end{proof}

\subsection*{Characteristic classes and the unitary version}
Define $\widetilde c_j(P)=c_j(u)$ and $\widetilde s_j(P)=s_j(u)$, with $s_j$ the Newton class and $s_0(u)=1$. For the compatible lift, $s_j(V)=k\,\widetilde s_j(P)$. On a base with $k$-torsion, \eqref{eq:two} specifies the normalized class before taking characteristic classes. In particular it supplies the information that division by $k$ in $H^{2j}(X;\Z)$ can lose. By the splitting principle the Newton formula for two classes $u,v$ of rank one is
\[
 s_j(uv)=\sum_{i=0}^j \binom jis_i(u)s_{j-i}(v).
\]
It applies to virtual bundles by additivity of the Newton classes and bilinearity of tensor product. Write $a_i=c_i(u), b_i=c_i(v)$. Newton's identities, first in the torsion-free universal ring $H^*(\BSU\times\BSU;\Z)$ and then by pullback, give
\begin{align*}
 c_2(uv)&=a_2+b_2, & c_3(uv)&=a_3+b_3,\\
 c_4(uv)&=a_4+b_4-5a_2 b_2, &
 c_5(uv)&=a_5+b_5-11(a_2b_3+a_3b_2).
\end{align*}

One can also retain an identified unitary lift of the projective core. Its parameter space is the homotopy pullback of $\Gr_{k,m}\to B\PU(k)\leftarrow\BU(k)$. Lemma~\ref{lem:pair} supplies a lift over $\Gr_{k,m}$, so the choices relative to the given projectivization are parametrized by $B\UU(1)=\CP^\infty$: tensor the distinguished vector factor by a line. Equivalently this trivializes the principal $B\UU(1)$-fibration obtained by pulling back the central extension. The distinguished factors are compatible with amplification up to isomorphism (on the homogeneous bundles they become direct sums), and tensoring two lifted factors tensors their extra lines. In the stable theory the decorated parameters give $\BSU_\otimes\times\CP^\infty\simeq\BU_\otimes$. Explicitly a rank-one virtual unit $z$ splits as $(z[\det z]^{-1},\det z)$, and the determinant is multiplicative on rank-one units. This explains the unitary extension in the appendix to \cite{old}.

\section{The core and its UHF stabilization}\label{sec:uhf}

Fix $k,\ell\ge2$ with $(k,\ell)=1$. Write $A_\ell=M_{\ell^\infty}$ and $D_{k,\ell}=M_k\otimes A_\ell$. Both are unital UHF algebras, the norm completions of the unions of $M_{\ell^N}$ and $M_{k\ell^N}$, respectively, under tensoring with identities. All tensor products of bundles are fibrewise. Let
\[
 G_{k,\ell}=\Aut(D_{k,\ell}),\qquad
 i:\PU(k)\longrightarrow G_{k,\ell},\quad
 i(\alpha)=\alpha\otimes\id_{A_\ell}.
\]
When a matrix bundle is described by its principal frame bundle, $Bi$ classifies unital extension of its structure group.

\begin{theorem}[Fixed-rank comparison]\label{thm:fibre}
With the matrix Grassmannians amplified in the second rank, there is a natural weak equivalence over $B\PU(k)$
\begin{equation}\label{eq:fibre}
 \Tel_{N\ge0}\Gr_{k,\ell^N}
 \ \simeq\ \hofib\!\left(B\PU(k)\xrightarrow{Bi} BG_{k,\ell}\right).
\end{equation}
Here $\Gr_{k,1}$ is a point, and the homotopy fibre is taken over the trivial bundle.
\end{theorem}
\begin{proof}
Put $E=\Emb(M_k,D_{k,\ell})$, $D=D_{k,\ell}$ and $P=\PU(k)$. On $E$ the right action of $P$ is precomposition. Restricting an automorphism of $D$ to the standard copy of $M_k$ gives a $P$-equivariant map $G_{k,\ell}\to E$ for the right action on the group by $i(P)$. This map is surjective. Indeed the images of a minimal projection under any two embeddings have equal trace $1/k$. Projections of equal trace in a UHF algebra are Murray--von Neumann equivalent (the trace is injective on $K_0$ and projections cancel in AF algebras). A partial isometry between the first diagonal projections extends to a unitary intertwining all matrix units: for systems $(e_{ij})$ and $(f_{ij})$ take $\sum_i f_{i1}v e_{1i}$ with $v^*v=e_{11}$ and $vv^*=f_{11}$. The resulting inner automorphism carries one embedding to the other.

Restriction has continuous local sections. For $(f_{ij})$ close to a fixed system $(e_{ij})$, the element $\sum_i f_{i1}e_{1i}$ is close to $1$, is invertible and intertwines the systems. Its polar unitary also intertwines them and is continuous in the entries. Compose the induced inner automorphism with an extension at the centre of the neighbourhood. The fibre at the standard inclusion consists of the automorphisms $\id_{M_k}\otimes\gamma$, $\gamma\in\Aut(A_\ell)$, as follows by restricting to the commutant. Restriction is therefore a locally trivial bundle. The group $\Aut(A_\ell)$ is contractible by \cite[Theorem~2.3]{DP}, so restriction is a weak equivalence.

At the finite stages $E_N=\Fr_{k,\ell^N}$ and at the limit use the orbit fibrations
\[
 \begin{aligned}
 \UU(\ell^N)&\longrightarrow\UU(k\ell^N)\longrightarrow E_N,\\
 \UU(A_\ell)&\longrightarrow\UU(D)\longrightarrow E .
 \end{aligned}
\]
The infinite-dimensional orbit map has local sections by the same polar formula. Both unitary-group telescopes map by weak equivalences to their limiting norm unitary groups. To see this, apply the conditional expectation onto a matrix stage to the entries of a continuous unitary family over a compact parameter space. The expectations converge uniformly on its compact image. At a sufficiently large stage take the polar part. Polar decomposition along the straight segment gives a homotopy to the original family. On a closed subset where the family already lies in that stage this homotopy is constant. Applied to spheres and relative disks, this proves the assertion including $\pi_0$. Passing to the direct limit in the exact homotopy sequences now gives $\Tel_N E_N\simeq E$ weakly.

Finally use the Borel construction for the right $P$-actions. The map $(E_N)_{hP}\to\Gr_{k,\ell^N}$ is a weak equivalence since the action is principal. The maps of these Borel constructions to $BP$ classify the frame bundles. Taking their telescope gives $E_{hP}$ up to weak equivalence, by the common Borel model and the fibre sequences over $BP$. The equivariant restriction equivalence identifies this with $(G_{k,\ell})_{hP}$. The latter over $BP$ is the principal $G_{k,\ell}$-bundle induced from $EP\to BP$ by $i$ (use inversion to exchange right and left multiplication). It is a model of the homotopy fibre of $Bi$ since it pulls back $EG_{k,\ell}\to BG_{k,\ell}$. These constructions also identify the map to $BP$ in \eqref{eq:fibre}.
\end{proof}

\begin{theorem}[Unital triviality and finite embeddings]\label{thm:criterion}
For an $M_k$-bundle $\Bun$ over $X$, the following conditions are equivalent.
\begin{enumerate}
\item $\Bun$ admits an embedding in $\underline M_{km}$ for some $m$ prime to $k$.
\item $\Bun\cong\End(V)$ for a rank-$k$ special-unitary bundle $V$ with $[V]=ku$ for some $u\in\units(X)$.
\item $\Bun$ embeds in $\underline M_{k\ell^N}$ for some $N\ge0$.
\item There is a unital algebra bundle isomorphism $\Bun\otimes A_\ell\cong X\times D_{k,\ell}$.
\end{enumerate}
If they hold, $\Bun$ embeds in $\underline M_{ks}$ for \emph{every} positive integer $s$ with $2s\ge\dim X$. In particular the assertion for powers of $\ell$ holds for any choice of $\ell\ge2$ prime to $k$.
\end{theorem}
\begin{proof}
Lemmas~\ref{lem:pair} and~\ref{lem:normalize} yield (1)$\Rightarrow$(2), including the constant type when $m=1$. Given (2), represent $su^{-1}$ by a vector bundle $W_s$ of rank $s$ when $2s\ge\dim X$. The product has class $[V\otimes W_s]=ks$, hence is trivial by cancellation. The action of $\End(V)$ on the first factor gives the embedding. This proves the stated bound and (2)$\Rightarrow$(3). Condition (3) implies (1) since $(k,\ell)=1$.

For $f:X\to BP$ classifying the core, condition (4) means that $Bi\circ f$ is nullhomotopic. This is equivalent to lifting $f$ up to homotopy into the homotopy fibre in \eqref{eq:fibre}. A map from finite $X$ to the telescope has bounded height and can be deformed into one stage by retracting a finite string of mapping cylinders. Here we can use a CW model for the homotopy fibre, since $X$ is a CW complex. By \eqref{eq:finite}, a lift into that stage is an embedded family whose core is isomorphic to $\Bun$. This proves (4)$\Leftrightarrow$(3).
\end{proof}

\begin{remark}
Condition (2) is existential in \emph{both} $V$ and $u$. For a fixed embedding the compatible pair specifies $u$, whereas distinct embeddings of an isomorphic core can give distinct stable classes. Section~\ref{sec:stages} gives their number on a Moore space. The equivalence \eqref{eq:fibre} also retains choices: a lift includes a nullhomotopy of the stabilized class. The vanishing criterion alone forgets this information.
\end{remark}

\subsection*{A finite test using two coprime core ranks}
The common-amplification test of \cite[Lemma~23]{old} has the following proof. Suppose that $\Bun_k\otimes\underline M_m\cong\mathcal C_m\otimes\underline M_k$ for bundles of the indicated matrix ranks and $(k,m)=1$. The two projective bundles have equal obstruction in $H^3(X;\Z)$ to a unitary lift. Indeed the scalar cocycle obtained from local unitary lifts of the transition functions is unchanged by amplification. The determinant shows that this class is killed by the corresponding matrix rank on either side, hence it is zero. Thus $\Bun_k=\End(V)$ and $\mathcal C_m=\End(W)$ for unitary bundles. The usual scalar obstruction can equivalently be obtained from $B\UU(1)\to B\UU(s)\to B\PU(s)$.

The amplified bundles agree projectively, so $V^{\oplus m}\cong W^{\oplus k}\otimes L$ for a line $L$. In particular $m c_1(V)-k c_1(W)=km c_1(L)$. Coprimality implies $c_1(V)\in kH^2(X;\Z)$ and $c_1(W)\in mH^2(X;\Z)$. Twist each factor by a line to make both determinants trivial. For the resulting comparison line, now denoted $L$, we have $L^{km}=\one$. Factor it as $L_kL_m$ with $L_k^k=L_m^m=\one$ by B\'ezout. Twisting $V$ by $L_k^{-1}$ and $W$ by $L_m$ gives actual special-unitary factors with $m[V]=k[W]$. For integers $a,b$ with $ak+bm=1$, the unit $u=a[V]+b[W]$ satisfies $[V]=ku$ and $[W]=mu$. Its two localized values $[V]/k$ and $[W]/m$ are independent of the compatible lifts, by the line argument in Lemma~\ref{lem:normalize}. This gives the unique common stable embedded class: one realizes it with either core by taking a sufficiently large coprime complementary rank and applying the construction in Theorem~\ref{thm:criterion}. Conversely, starting with $\Bun_k$ and its class $u$, for any $2m\ge\dim X$ represent $mu$ by $W$. Stable cancellation gives $V^{\oplus m}\cong W^{\oplus k}$ and hence the common amplification. This proves the finite test and the uniqueness statement for the pair of cores, without uniqueness of a lift of one core.

\subsection*{What the two limits remember}
The frame computation also gives the homotopy groups of the unital automorphism group. Write $R_\ell=\Z[1/\ell]$. The stable odd homotopy groups of the two norm unitary groups in the proof of Theorem~\ref{thm:fibre} are $R_\ell$, and the positive even groups vanish. In rank-one matrix coordinates, normalized by $\ell^N$ in both towers, $u\mapsto1_k\otimes u$ acts as multiplication by $k$. The groups are path connected, so the exact sequence and restriction equivalence give
\[
 \pi_0G_{k,\ell}=0,\qquad
 \pi_{2j-1}G_{k,\ell}=R_\ell/kR_\ell\cong\Z/k,\qquad
 \pi_{2j}G_{k,\ell}=0\quad(j\ge1).
\]
For two complementary parameters $\ell,h$ prime to $k$, tensoring with $\id_{A_h}$ maps $G_{k,\ell}$ to $G_{k,\ell h}$ using the interleaved tensor factors. On the frame computation it induces $R_\ell/kR_\ell\to R_{\ell h}/kR_{\ell h}$: the finite-stage multiplicity $h^N$ is cancelled by the new denominator. This is an isomorphism in every degree. Interchanging $h$ and $\ell$ gives a zigzag of weak equivalences of groups compatible with the inclusions of $\PU(k)$. Thus the resulting fibre models over $B\PU(k)$ are independent of the complementary parameter up to weak equivalence. The earlier frame-space independence appears in \cite[Lemma~28]{brauerpre}.

There is a separate finite-matrix localization in \cite[Section~3.4]{old}. For $p\ge2$ let $T_p=\Tel_n B\PU(p^n)$, with amplification by $p$. It classifies matrix bundles of $p$-power rank up to amplification by trivial $p$-power matrices on finite bases. This follows by factoring maps and homotopies through finite telescope segments. The group $\pi_2(T_p)$ is $\Z[1/p]/\Z$ by the direct limit of the groups $\Z/p^n$, while $\Aut(A_p)$ is contractible. In degrees $2j\ge4$ one obtains $\pi_{2j}(T_p)=\Z[1/p]$ and the odd groups vanish, directly by Bott periodicity. With unrestricted matrix ranks the corresponding groups are $\Q/\Z$ in degree two and $\Q$ in higher even degrees, as seen along the factorial tower.

To describe the forgetful map to $[X,T_p]$, realize $u$ using a first factor of sufficiently large rank $p^n$ and then retain its endomorphism algebra. This gives the amplified core with vector class $p^n u$. Equal normalized classes give the same amplified core by cancellation for the distinguished vector representatives. The image is the trivial class precisely when $u=1$ in $K^0(X)[1/p]$. Indeed if the amplified endomorphism bundle is trivial then its vector factor is a sum of copies of a line. Its trivial determinant makes this line of $p$-power order, so it becomes trivial in $K^0(X)[1/p]$. Conversely if $p^a(u-1)=0$, sufficient amplification gives a stably trivial vector representative in the cancellation range. The kernel thus consists of the classes for which $u-1$ is killed by a power of $p$. For the unrestricted tower the corresponding condition is additive torsion of $u-1$. In multiplicative language these are the same respective torsion subgroups: the rational nilpotent augmentation ideal has torsion-free group of units $1+x$ by the logarithm, and a finite additive torsion ideal generated by $x$ has unit group filtered by the additive quotients of its powers. For the $p$-power version, if $p^a x=0$ then the ideal generated by $x$ is a finite group killed by $p^a$, and its unit filtration gives the same restriction on prime divisors. Conversely, if $u^{p^a}=1$, the factor $1+u+\cdots+u^{p^a-1}$ is a unit after inverting $p$, by rank and nilpotence, hence $u=1$ there. In the rank-preserving variant in \cite[Remark~25]{old}, the values at a point are initially the monoids of allowed positive ranks. Group completion gives $\Q_{>0}^{\times}$ for all ranks and $p^\Z$ for $p$-powers.

\section{Least finite stages and torsion}\label{sec:stages}

Let $\Bun=\End(V)$ be a fixed split core of matrix rank $k$. Define its \emph{least $\ell$-stage} to be the least $N\ge0$ for which it embeds in $\underline M_{k\ell^N}$, with value $\infty$ if the set is empty. The core rank stays fixed in this definition. The proof of Theorem~\ref{thm:criterion} bounds the stage whenever a coprime embedding exists. To decide sharpness one must take account of the complementary rank.

\begin{lemma}[Multiplicity bundle]\label{lem:mult}
An embedding $\End(V)\to\underline M_{km}$ exists precisely when there is a rank-$m$ bundle $W$ with $V\otimes W\cong\one^{km}$. When $\dim X\le2km$, this is equivalent to the existence of a rank-$m$ vector representative of some $w\in K^0(X)$ with $[V]w=km$.
\end{lemma}
\begin{proof}
Regard $\one^{km}$ as a left $\End(V)$-module via the embedding. The intertwiner spaces
\[
 W_x=\operatorname{Hom}_{\End(V_x)}(V_x,\C^{km})
\]
have constant dimension $m$ by the representation theory of $M_k$. They form a vector bundle: in a chart they are the kernel of a continuous bundle map imposing the intertwining relations on finitely many matrix units, of constant rank. Evaluation $V\otimes W\to\one^{km}$ is an isomorphism on every fibre. Conversely, a unitary trivialization of the tensor product realizes the required action on the first factor. A vector bundle trivialization can always be made unitary by polar decomposition. In the stated dimension range zero reduced class of $V\otimes W$ implies actual triviality.
\end{proof}

Twisting $V$ by a line $L$ carries the multiplicity solutions bijectively to those for the new lift by $W\mapsto W\otimes L^{-1}$. In particular the least stage depends only on the algebra bundle. Condition (2) of Theorem~\ref{thm:criterion} is convenient for eventual existence, while Lemma~\ref{lem:mult} applies to a prescribed multiplicity with no coprimality assumption.

\begin{proposition}[All stable classes with a prescribed core]\label{prop:ambiguity}
Suppose $\Bun$ is an embeddable $M_k$-bundle and one embedding gives the normalized unit $u_0$. Put
\[
 P_k(X)=\{1+a:a\in K^0(X),\ ka=0\},\qquad
 d(y)=y[\det y]^{-1}\quad (y\in P_k(X)).
\]
Then the set of stable embedded classes admitting a representative with abstract core $\Bun$ corresponds exactly to the coset
\[
                         u_0\,d(P_k(X))\ \subset\ \units(X).
\]
In particular the subgroup giving the ambiguity depends on $k$ and $X$ alone.
\end{proposition}
\begin{proof}
Rank considerations and nilpotence show that $P_k(X)$ consists of units of rank one. It is a group since $y\in P_k(X)$ means $k(y-1)=0$, a condition preserved by products and inverses. The determinant is multiplicative for tensor products of virtual rank-one classes, so $d$ is a homomorphism. If $y=1+a$, then $(\det y)^k=\one$ by additivity of the first Chern class and $ka=0$.

Let $V_0$ be a compatible special-unitary vector lift with $[V_0]=ku_0$. Any other compatible lift is $V_0\otimes L$ for a line with $L^k=\one$. If its unit is $u$, we have $ku=ku_0[L]$. Thus $y=u u_0^{-1}[L]^{-1}$ lies in $P_k(X)$ and $\det y=L^{-1}$, giving $u=u_0d(y)$.

Conversely let $y\in P_k(X)$ and put $L=(\det y)^{-1}$. The bundle $V_0\otimes L$ has trivial determinant and endomorphism algebra $\Bun$. Its class equals $k u_0d(y)$ since $ky=k$. For $m$ coprime to $k$ with $2m\ge\dim X$, choose a rank-$m$ bundle of class $m(u_0d(y))^{-1}$. Its tensor product with $V_0\otimes L$ is trivial by cancellation. This embedding has normalized class $u_0d(y)$: the prescribed vector factors have the two local values in Lemma~\ref{lem:normalize}, unchanged on passing to their distinguished projective lifts.
\end{proof}

\begin{theorem}[All multiplicities on an even sphere]\label{thm:sphere}
Let $X=S^{2j}$ with $2\le j\le k$. Fix a Bott generator $\beta\in\widetilde K^0(S^{2j})$ and write $[V]=k+t\beta$. If $t\ne0$, then
\begin{equation}\label{eq:sphere}
  \End(V)\hookrightarrow\underline M_{km}
   \qquad\Longleftrightarrow\qquad
   m\ge j\ \text{ and }\ k\mid mt.
\end{equation}
If $t=0$, every positive multiplicity works. For $(k,\ell)=1$, the least $\ell$-stage is $0$ for $t=0$, is $\infty$ for $k\nmid t$, and otherwise is the least $N$ with $\ell^N\ge j$.
\end{theorem}
\begin{proof}
The sphere is a suspension and $\beta^2=0$. Write $[W]=m+s\beta$ for a multiplicity bundle. Then $[V][W]=km$ gives $mt+ks=0$. If $m<j$, all possible positive Chern classes of $W$ vanish: their groups below degree $2j$ are zero, and $c_j(W)=0$ by rank. Consequently its reduced Chern character vanishes, hence $s=0$ since $\widetilde K^0(S^{2j})$ is torsion-free. This would give $t=0$. Conversely, if $m\ge j$ and $k\mid mt$, the integral class $m-(mt/k)\beta$ is realizable in rank $m$ by stability. The tensor product with $V$ has trivial stable class in the cancellation range, so Lemma~\ref{lem:mult} gives the embedding. For $t=0$ the bundle $V$ is already trivial since $k\ge j$. Coprimality gives the assertions for $m=\ell^N$.
\end{proof}

The integer $t$ in \eqref{eq:sphere} is a stable clutching number. With a choice of orientation, $c_j(V)[S^{2j}]=(-1)^{j-1}(j-1)!\,t$ by the Chern character of the Bott class. Divisibility of this Chern number by $k$ alone can therefore lose information. For example, on $S^8$ a rank-six bundle of Bott integer $1$ has top Chern number $\pm6$, but it has no coprime embedding.

\begin{proposition}[A four-dimensional Moore space]\label{prop:moore}
Let $X_q=S^3\cup_q e^4$, $q\ge2$, and $k\ge2$. In rank $k$ write $[V]=k+t\gamma$, where $\gamma$ generates $\widetilde K^0(X_q)=\Z/q$ and $t\in\Z/q$. Put $g=(k,q)$. For $m\ge2$ the algebra $\End(V)$ embeds in $\underline M_{km}$ precisely when $g\mid mt$ in the sense of integer representatives. Multiplicity one works precisely when $t=0$. For every $\ell\ge2$ prime to $k$ the least $\ell$-stage is therefore
\[
 \begin{cases}
 0,& t=0,\\
 1,& t\ne0,\quad g\mid t,\\
 \infty,& g\nmid t.
 \end{cases}
\]
If $g\mid t$, exactly $g$ stable embedded classes have representatives with this given core.
\end{proposition}
\begin{proof}
The cofibration for the degree-$q$ map of $S^3$ and Bott periodicity give $\widetilde K^0(X_q)=\Z/q$ and $K^1(X_q)=0$. The space is a suspension, so $\gamma^2=0$. Its $H^2$ vanishes, hence it has no nontrivial line bundles. A rank-$m$ solution for $m\ge2$ is therefore equivalent by stability and Lemma~\ref{lem:mult} to a residue $s$ with $mt+ks=0$ in $\Z/q$. The image of multiplication by $k$ consists of multiples of $g$. For $m=1$ the multiplicity line must be trivial, and stable cancellation in rank $k$ gives the assertion. Since $(\ell,g)=1$, the displayed stages follow.

All vector lifts of the given projective core are isomorphic, since $H^2(X_q;\Z)=0$. A possible normalized unit has the form $1+x\gamma$ and must satisfy $kx=t$ in $\Z/q$. Conversely, for each such $x$ choose any complementary rank $m\ge2$ prime to $k$. Take a representative of $m(1-x\gamma)$ and use stable cancellation to trivialize its tensor product with $V$. The pair obtained has the given unit by \eqref{eq:two}. There are precisely $g$ solutions, distinguished by Theorem~\ref{thm:stable}.
\end{proof}

For instance, the trivial $M_3$-bundle on $X_3$ occurs in three distinct stable classes, with units $1,1+\gamma,1+2\gamma$. Here the vector lift itself is unique. Here $d$ is the identity on $P_3(X_3)$ in Proposition~\ref{prop:ambiguity}. The information lost on forgetting the embedding lies in the possible complements. Example~\ref{ex:RP} exhibits a separate effect: ambiguity of the vector lift.

\section{From the stable unit to a higher twist}\label{sec:twists}

There are two operator-algebra constructions associated to an embedding. The core gives the unital bundle $\Bun\otimes A_\ell$ on $X$, whose triviality is characterized by Theorem~\ref{thm:criterion}. Its full embedded class supplies a unit $u$ and hence a suspension construction over $\Sigma X$. We spell out the latter using complex $K$-theory localized at $\ell$.

Let $\Kc$ be the compact operators on an infinite-dimensional separable Hilbert space. The group isomorphism of Dadarlat and Pennig \cite[Theorem~2.22]{DP} sends a homotopy class of maps $X\to\Aut(A_\ell\otimes\Kc)$ to its evaluation on the projection $1\otimes e$, with $e$ of rank one. Its target is the positive units of
\[
 K_0(C(X)\otimes A_\ell)=K^0(X)[1/\ell].
\]
The equality follows by continuity from the matrix stages in rank-one coordinates. The image of our integral unit $u$ and its inverse are positive: for large $N$ each of $\ell^Nu$ and $\ell^Nu^{-1}$ is represented by a vector bundle, and the denominators are removed by passing to $A_\ell$. Thus $u$ gives a clutching map $\alpha_u:X\to\Aut(A_\ell\otimes\Kc)$, unique up to homotopy. Gluing the trivial bundles on the two cones defines a twist $\tau_u$ on $\Sigma X$. We use the identification with twists of $KU[1/\ell]$ from \cite{unitspectra,Pennig}. Reversing which cone is glued first inverts the unit.

\begin{proposition}\label{prop:twist}
For a coprime embedding on $S^{2j}$ as in Theorem~\ref{thm:sphere}, write $t=kh$. Then $u=1+h\beta$. If $h\ne0$, the induced twist on $S^{2j+1}$ has groups
\[
 KU[1/\ell]^0(S^{2j+1};\tau_u)=0,\qquad
 KU[1/\ell]^1(S^{2j+1};\tau_u)\cong R_\ell/h R_\ell .
\]
\end{proposition}
\begin{proof}
Equation~\eqref{eq:two} determines $u=1+h\beta$ because $K^0(S^{2j})$ is torsion-free. Over the cones the coefficient groups in degree zero are $R_\ell$ and their odd groups vanish. On the equator they are $R_\ell\oplus R_\ell\beta$ and zero, respectively. With one choice of gluing convention the Mayer--Vietoris difference map is
\[
 R_\ell^2\longrightarrow R_\ell\oplus R_\ell\beta,\qquad
 (a,b)\longmapsto(a-b)-bh\beta .
\]
Indeed the change of trivialization sends the constant rank-one class on the equator to the evaluation class $u$. For $h\ne0$ the kernel is zero and the cokernel is $R_\ell/hR_\ell$. This is the suspension computation of \cite[Section~4]{Pennig} in the present coordinates.
\end{proof}

\begin{example}\label{ex:S6}
Choose $j=k=3$, $\ell=2$ and $[V]=3+9\beta$ on $S^6$. In rank two a multiplicity bundle would have vanishing reduced class, contradicting $2\cdot9+3s=0$. In rank four take $[W]=4-12\beta$. Then $V\otimes W$ is trivial. Thus stage $M_6$ fails, stage $M_{12}$ works, and the least exponent in the dyadic tower is $2$. The full pair has unit $1+3\beta$, and $\End(V)\otimes M_{2^\infty}$ is a trivial unital bundle on $S^6$. Clutching by the localized unit gives a twist on $S^7$ with even group zero and odd group $\Z/3$. Its nonzero clutching class belongs to $\pi_6\Aut(M_{2^\infty}\otimes\Kc)\cong\Z[1/2]$, while $H^3(S^7;\Z)=0$. Thus the example realizes a higher twist and specifies the first finite matrix stage for this core.
\end{example}

\paragraph{Further questions.}
For a prescribed core on a general finite CW complex, Lemma~\ref{lem:mult} reduces the least-multiplicity problem, in the ambient cancellation range, to the geometric dimensions of the solutions of $[V]w=km$. The Moore space computation uses that every rank at least two lies in the stable range and that every reduced product vanishes. On a higher-dimensional lens space, reduced products and possible low-rank representatives both enter this equation. Proposition~\ref{prop:ambiguity} describes all stable units over a fixed core. A finer question is which members of that coset can already be realized at a \emph{specified} complementary rank, and how distinct embeddings at that rank can be distinguished before stabilization. These questions concern finite representatives and remain meaningful after the UHF obstruction has vanished.

\paragraph{Acknowledgments.}
I am grateful to the Max Planck Institute for Mathematics in Bonn for its support and hospitality during part of the original work. I thank E.\,V. Troitsky for his continued interest and support, and V.\,M. Manuilov and A.\,S. Mishchenko for discussions. The original research was supported in part by RFFI grant 02-01-00572.

\paragraph{Use of AI tools.}
OpenAI's ChatGPT assisted with reviewing the earlier arguments, deriving and checking examples, and preparing this revised manuscript.

\end{document}